\documentclass[12pt]
{amsart}
\usepackage[cp1251]{inputenc}
\usepackage{amsmath,amssymb,amsthm,comment,graphicx}
\usepackage[unicode, pdftex]{hyperref}
\usepackage[russian]{babel}

\newcommand{\R}{\mathbb{R}}
\newcommand{\N}{\mathbb{N}}
\newcommand{\Z}{\mathbb{Z}}
\newcommand{\Q}{\mathbb{Q}}
\newcommand{\K}{\mathbb{K}}
\newcommand{\eps}{\varepsilon}

\newtheorem{theorem}{Теорема}
\newtheorem{lemma}{Лемма}
\newtheorem{prop}{Предложение}
\newtheorem*{question}{Вопрос}

\theoremstyle{remark}
\newtheorem*{remark}{Замечание}

\theoremstyle{definition}

\newtheorem{problem}{Задача}

\begin{document}

\title{Разрезание треугольника на равные треугольники}
\author{А. Рябичев}
\thanks{}
\date{}
\maketitle

\begin{abstract}
В статье доказывается, что {\it почти любой} треугольник
можно разрезать лишь на $n^2$ равных между собой треугольников,
причём для каждого $n$ такое разбиение единственно.
Для решения этой <<школьний>> задачи мы используем чуть более продвинутую технику,
чем просто школьную геометрию, --- некоторые приёмы из линейной алгебры, анализа и теории меры.
\end{abstract}

\section{Введение}

\begin{figure}[h]
\includegraphics[width=0.3\textwidth]{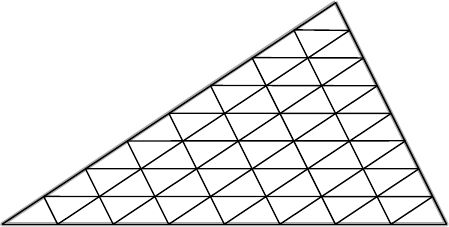}
\caption{Разрезание треугольника на $n^2$ равных треугольников.}
\label{fig:n2-cut}
\end{figure}

Любой треугольник можно разрезать на $n^2$ равных между собой треугольников для любого $n\in\N$.
Схема такого разрезания
показана на рис.~\ref{fig:n2-cut}.
Будем называть такие разбиения {\it стандартными}.

Некоторые треугольники также допускают другие разрезания
на $m$ равных между собой треугольников, в том числе для чисел $m$, не являющихся полными квадратами.
Мы предлагаем читателю рассмотреть самые простые примеры, приведённые ниже, самостоятельно.


\begin{problem}
Докажите, что если треугольник $\Delta$ можно разбить на 2 равных между собой треугольника, то $\Delta$ равнобедренный.
\end{problem}

\begin{problem}
{\bf (а)} Докажите, что если треугольник $\Delta$ можно разбить на 3 равных между собой треугольника, то $\Delta$ имеет угол $\pi/3$.
{\bf (б)} Докажите, что таких треугольников всего два (с точностью до подобия).
\end{problem}

\begin{problem}
Постройте треугольник, который можно разбить на 4 равных между собой треугольников способом, отличным от стандартного.
\end{problem}

\begin{problem}
Существует ли треугольник, который можно разбить на 5 равных между собой треугольников?
\end{problem}

При большом $m$ анализ способов разбить треугольник на $m$ равных между собой треугольников довольно сложен.
В нашей статье мы поговорим о следующей задаче.

\begin{question}
Существует ли треугольник, который можно разрезать лишь на $n^2$ равных между собой треугольников?
\end{question}

В книге \cite{S90} А.\,Сойфер даёт положительный ответ на этот вопрос,
презентуя ряд изящных и весьма полезных приёмов из линейной алгебры,
так что приведённое им доказательство вполне доступно заинтересованному старшекласснику.

Мы дадим немного другое доказательство существования такого треугольника.
При этом мы выходим за рамки школьной геометрии и используем чуть более продвинутые приёмы из анализа и теории меры.
Как нам кажется, эта техника довольно естественна для решения нашей задачи.
В частности, она легко приводит нас к более общему результату, чем просто теорема существования:

\begin{theorem}\label{th:cut}
Почти любой треугольник допускает разбиение лишь на $n^2$ равных между собой треугольников для всех натуральных $n$.
\end{theorem}

Более тщательный анализ разбиений даёт следующее замечательное обобщение:

\begin{theorem}\label{th:unique}
Почти для любого треугольника все разбиения на равные между собой треугольники
могут быть устроены только как на Рис.~\ref{fig:n2-cut}.
\end{theorem}

В формулировках теорем \ref{th:cut} и \ref{th:unique} требуется уточнить,
что означают слова {\it почти любой}. 
В следующих двух разделах мы даём необходимые определения.

\section{Множества меры ноль}

Говорят, что множество $A\subset\R^k$ имеет {\it меру ноль},
если для любого $\eps>0$ множество $A$ можно покрыть
не более чем счётным числом параллелепипедов
со сторонами, параллельными координатным осям, и имеющими суммарный объём меньше $\eps$.
В этих терминах, фраза <<{\it почти любая точка~$\R^k$}>> имеет буквальное значение
<<{\it любая точка $\R^k$, не принадлежащая наперёд заданному множеству меры ноль}>>.

Легко показать, что объединение не более чем счётного числа множеств меры ноль будет иметь меру ноль.
Также несложно убедиться, используя компактность, что никакой шар в $\R^k$ не является множеством меры ноль.
(См., например, \cite[стр.~337]{zorich1}.)

Мы будем использовать два следующих факта.
Их можно считать упражнениями по введению в теорию меры.

\begin{prop}\label{pf:plane-zero}
Объединение счётного числа прямых в $\R^2$ имеет меру ноль.
\end{prop}

Пусть даны открытые подмножества $U,V\subset\R^2$.
Мы называем отображение $f:U\to V$ {\it гладким},
если его координатные функции непрерывно-дифференцируемы.
Следующее утверждение нетрудно выводится из локальной липшицевости гладкого отображения.

\begin{prop}\label{pf:smooth-image}
Для любого множества меры ноль $A\subset\R^2$
и любого гладкого отображения $f:U\to V$
множество $f(A\cap U)$ имеет меру ноль.
\end{prop}

\section{Пространство треугольников}

Каждому треугольнику $\Delta$ сопоставим упорядоченную тройку длин его сторон $(a,b,c)$.
Ясно, что для нашей задачи нам достаточно рассмотреть лишь те треугольники, у которых одна из сторон равна единице.
Поэтому будем считать, что $c=1$.

Таким образом, множество всех треугольников (с упорядоченными вершинами, с точностью до подобия)
отождествляется с подмножеством $M\subset\R^2$, заданным неравенствами $x+1>y$, \ $y+1>x$ \ и \ $x+y>1$.

Далее мы построим множество $\Sigma\subset M$ меры ноль, такое что
для любого $\Delta\in M\setminus\Sigma$ выполнено утверждение теорем \ref{th:cut} и \ref{th:unique}.
Поскольку дополнение $M\setminus\Sigma$ непусто, теорема А.\,Сойфера сводится к теореме~\ref{th:cut}.

\section{Несоизмеримость углов}

Пусть $N\subset\R^2$ --- подмножество, заданное уравнениями $x>0$, $y>0$ и $x+y<\pi$.
Для $(\alpha,\beta)\in N$ существует единственный треугольник $\Delta\in M$
с упорядоченным набором углов $(\alpha,\beta,\pi-\alpha-\beta)$.
Это задаёт отображение $f:N\to M$.

\begin{prop}
Отображение $f$ гладко.
\end{prop}

\begin{proof}
Возьмём произвольный треугольник $ABC$, такой что $AB=1$.
Тогда по теореме синусов $AC=\frac{\sin\angle B}{\sin\angle C}$ и $BC=\frac{\sin\angle A}{\sin\angle C}$.
Отношение гладких функций гладко, поэтому зависимость сторон треугольника $ABC$ от его углов является гладкой, что и требовалось.
\end{proof}

Для каждой тройки рациональных чисел $(q_1,q_2,q_3)\ne(0,0,0)$ уравнение
$$q_1x+q_2y+q_3(\pi-x-y)=0$$ задаёт прямую, либо пустое множество.
Пусть $\Sigma_1\subset\R^2$ --- объединение всех прямых такого вида.
Тогда по предложению \ref{pf:plane-zero} множество $\Sigma_1$ имеет меру ноль,
а по предложению \ref{pf:smooth-image} множество $f(\Sigma_1\cap N)$ имеет меру ноль.

Про треугольники из $M\setminus f(\Sigma_1\cap N)$ говорят, что {\it их углы не соизмеримы}.

\begin{remark}
Если $\alpha,\beta,\gamma$ --- набор углов некоторого треугольника,
такой что $(\alpha,\beta)\in\Sigma_1$, то тогда также $(\beta,\alpha)\in\Sigma_1$, $(\alpha,\gamma)\in\Sigma_1$, и т.\,д.\
(то есть не важно, какую пару углов брать).
\end{remark}

\begin{lemma}\label{l:cong}
Пусть $\Delta\in M\setminus f(\Sigma_1\cap N)$.
Тогда, если $\Delta$ разбит на треугольники, равные~$\Delta'$,
то $\Delta$ и $\Delta'$ подобны.
\end{lemma}

\begin{proof}
Предположим обратное, $\Delta$ разбит на треугольники, равные $\Delta'$, и при этом $\Delta'$ не подобен $\Delta$.
Пусть $\alpha<\beta<\gamma$ --- углы $\Delta$,
а $\alpha'\le\beta'\le\gamma'$ --- углы $\Delta'$.

Заметим, что обязательно выполнено хотя бы одно из неравенств:
$\alpha'>\alpha$, $\beta'>\beta$ или $\gamma'>\gamma$.
В противном случае либо $\alpha'+\beta'+\gamma'<\alpha+\beta+\gamma$, либо $\Delta'$ и $\Delta$ окажутся подобны.

Если $\alpha'>\alpha$, то внутри угла $\alpha$ не поместится ни один из углов $\Delta'$, что противоречит наличию разбиения.

Если $\beta'>\beta$, то, из существования разбиения, внутри $\alpha$ и $\beta$ помещается
целое число углов, равных $\alpha'$, откуда $(\alpha,\beta)\in\Sigma_1$.

Наконец, пусть $\gamma'>\gamma$.
Тогда каждый угол $\Delta$ может быть составлен из углов $\alpha',\beta'$.
Другими словами, для некоторых $u_1,u_2,v_1,v_2,w_1,w_2\in\Z_{\ge0}$ мы имеем
$$
\begin{cases}
\alpha=u_1\alpha'+u_2\beta'\\
\beta=v_1\alpha'+v_2\beta'\\
\gamma=w_1\alpha'+w_2\beta'
\end{cases}
$$
Но вектора $u=(u_1,u_2)$, $v=(v_1,v_2)$ и $w=(w_1,w_2)$ в $\Q^2$ линейно зависимы.
Это значит, что $q_1u+q_2v+q_3w=0$ для некоторой тройки рациональных чисел $(q_1,q_2,q_3)\ne(0,0,0)$.
Отсюда следует равенство $q_1\alpha+q_2\beta+q_3\gamma=0$, что противоречит предположению $\Delta\notin f(\Sigma_1\cap N)$.
\end{proof}

\section{Несоизмеримость сторон}

Пусть $\K\subset\R$ --- наименьшее подполе, содержащее квадратные корни всех натуральных чисел.
Опишем его явно: все элементы $\K$ суть дроби вида
$\frac{\pm\sqrt m_1\pm\ldots\pm\sqrt m_k}{\pm\sqrt n_1\pm\ldots\pm\sqrt n_l}$ с ненулевыми знаменателями,
где $m_1,\ldots,m_k,n_1\ldots,n_l\in\N$.
Очевидно, $\K\supset\Q$ и $\K$ счётно.

Для каждой тройки $q_1,q_2,q_3$ элементов $\K$, такой что $(q_1,q_2)\ne(0,0)$,
возьмём прямую в $\R^2$, заданную уравнением $q_1x+q_2y+q_3=0$.
Пусть $\Sigma_2\subset\R^2$ --- объединение всех прямых такого вида.

Про треугольники из $M\setminus\Sigma_2$ говорят, что {\it их стороны не соизмеримы над $\K$}.
По предложению \ref{pf:plane-zero} множество $\Sigma_2$ имеет меру ноль.

\begin{lemma}\label{l:n2-cong}
Пусть $\Delta\in M\setminus(f(\Sigma_1\cap N)\cup\Sigma_2)$.
Тогда, если $\Delta$ разбит на треугольники, равные~$\Delta'$,
то $\Delta$ и $\Delta'$ подобны с коэффициентом $\frac1{n^2}$ для некоторого $n\in\N$.
\end{lemma}

\begin{proof}
Пусть $\Delta$ разбит на $m$ копий треугольника $\Delta'$.
По лемме~\ref{l:cong} треугольники $\Delta$ и $\Delta'$ подобны.
Из соотношения площадей, коэффициент подобия равен $\frac1{\sqrt m}$.

Обозначим стороны $\Delta$ через $a,b,c$, причём $c=1$.
Тогда $\Delta'$ имеет стороны $\frac1{\sqrt m}a$, $\frac1{\sqrt m}b$, $\frac1{\sqrt m}c$.
Рассмотрим замощение стороны $c$.
Для некоторых $n_1,n_2,n_3\in\Z_{\ge0}$ мы имеем
$$
c=n_1\textstyle\frac1{\sqrt m}a+n_2\frac1{\sqrt m}b+n_3\frac1{\sqrt m}c,
$$
откуда $\frac{n_1}{\sqrt m}a+\frac{n_2}{\sqrt m}b+\frac{n_3-\sqrt m}{\sqrt m}=0$.
Из предположения $\Delta\notin\Sigma_2$, оба числа $n_1,n_2$ равны нулю.
Поэтому, раз равенство верно, $n_3=\sqrt m$, то есть $m=n_3^2$.
\end{proof}

\begin{proof}[Доказательство теоремы~\ref{th:cut}]
Нам достаточно положить $\Sigma=(f(\Sigma_1\cap N)\cup\Sigma_2)\cap M$.
Действительно, $\Sigma$ имеет меру ноль и, согласно лемме~\ref{l:n2-cong},
любой треугольник $\Delta\in M\setminus\Sigma$ можно разбить лишь на $n^2$ равных между собой треугольников.
\end{proof}


\section{Области со стандартным разбиением}

Начнём техническую подготовку к доказательству теоремы~\ref{th:unique}.
Пусть дано разбиение некоторого подмножества плоскости на равные между собой треугольники.
Будем называть разбиение {\it стандартным},
если все треугольники лежат в некоторой решётке как на рисунке~\ref{fig:n2-cut}.

Зафиксируем треугольник $\Delta\in M$.
Возьмём ограниченное множество $U\subset\R^2$, допускающее стандартное разбиение на треугольники, равные $\Delta$.

\begin{figure}[h]
\includegraphics[width=0.35\textwidth]{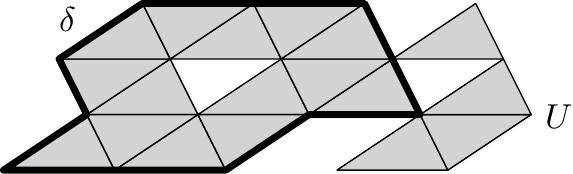}
\caption{Компонента границы $U$.}
\label{fig:delta}
\end{figure}

Пусть $\delta$ --- компонента границы $U$.
Более точно, $\delta$ является замкнутой ломаной,
по одну сторону к $\delta$ по всей длине примыкает $U$,
а по другую сторону вблизи $\delta$ лежит незаполненная часть плоскости
и, возможно, некоторые треугольники, которые соприкасаются с $\delta$ лишь вершиной, см.\ рис.~\ref{fig:delta}.
Также предположим, что одна из компонент $\R^2\setminus U$, примыкающая к $\delta$ хотя бы по одному ребру, неограничена.

Будем считать вершинами $\delta$ только те точки, где она делает излом.
Углы $\delta$ будем мерить так, что вблизи вершины внутренняя часть угла покрыта областью $U$.
Таким образом, углы $\delta$ принимают значения в $(0;\pi)\cup(\pi;2\pi)$.
Выберем направление обхода $\delta$.

\begin{lemma}\label{l:2or1sharp}
У $\delta$ найдутся либо два угла $<\pi$, идущих подряд, либо угол $>\pi$, оба соседних с которым $<\pi$.
\end{lemma}

\begin{proof}
Заметим, что для рёбер $\delta$ есть всего шесть вариантов направления.
Отметим на единичной окружности 6 точек, соответствующих этим направлениям.
Мы получили своего рода циферблат с шестью <<делениями>>.

Если мы проходим угол $<\pi$, следующее ребро $\delta$ поворачивается относительно предыдущего на $1$ или $2$ деления,
а если мы проходим угол $>\pi$, то следующее ребро $\delta$ поворачивается относительно предыдущего на $-1$ или $-2$ деления.

\begin{figure}[h]
\includegraphics[width=0.35\textwidth]{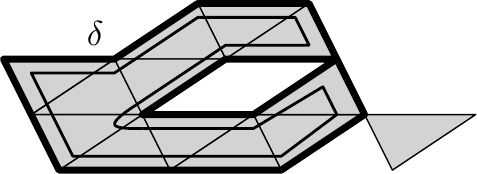}
\caption{Компонента границы $U$ и идущая вдоль неё кривая внутри $U$.}
\label{fig:delta2}
\end{figure}

Вообще говоря, $\delta$ может иметь самопересечения.
Но, если взять кривую, идущую вдоль $\delta$ внутри $U$, она будет несамопересекающейся, см.\ рис.~\ref{fig:delta2}.
Поэтому, суммируя количество делений
(на которое следующее ребро повернулось относительно предыдущего)
по всем углам $\delta$, за один полный круг мы получим ровно $6$.

Однако, если после каждого угла $<\pi$ идут хотя бы два угла $>\pi$,
суммарный вклад этих трёх углов в число делений будет $\le0$.
Поэтому утверждение леммы должно выполняться.
\end{proof}

\section{Доказательство теоремы~\ref{th:unique}}

Возьмём ограниченное множество $U\subset\R^2$ и треугольник $\Delta\in M\setminus\Sigma$.
Предположим, $U$ допускает стандартное разбиение на треугольники, равные $\Delta$.
Определим ломаную $\delta$ как в предыдущем разделе.

Предположим, $U$ также имеет разбиение на треугольники, равные $\Delta$, не совпадающее со стандартным разбиением.
Обозначим это разбиение~$P$.

\begin{prop}\label{prop:2sharp}
Пусть $I$ --- ребро $\delta$, оба угла при котором $<\pi$.
Тогда хотя бы один из треугольников разбиения $P$, примыкающий к $I$, совпадает с треугольником из стандартного разбиения.
\end{prop}

\begin{proof}
Пусть $a,b,c$ --- длины сторон $\Delta$, пусть для определённости $|I|=ma$.
Поскольку $a$, $b$, и $c$ не соизмеримы, все треугольники разбиения $P$ примыкают к $I$ стороной~$a$.

Пусть $\alpha,\beta,\gamma$ --- углы $\Delta$, лежащие против сторон $a,b,c$, соответственно.
Возьмём вершину $B$ ребра $I$.
Пусть треугольник стандартного разбиения, примыкающий к $I$ и содержащий $B$, имеет там угол $\beta$.
Тогда угол ломаной $\delta$ в вершине $B$ равен либо $\beta$, либо $\beta+\alpha$
(в~зависимости от того, 1 или 2 треугольника стандартного разбиения он содержит).

Из несоизмеримости $\alpha$, $\beta$ и $\gamma$, угол ломаной $\delta$ в вершине $B$
может быть представлен в виде комбинации $\alpha$, $\beta$ и $\gamma$ с рациональными коэффициентами единственным способом.
Поэтому треугольник разбиения $P$, примыкающий к $I$ и содержащий $B$, не может иметь там угол $\gamma$.
Следовательно, этот треугольник расположен как в стандартном разбиении.

(Итерируя этот аргумент, можно показать, что
вообще все треугольники $P$, примыкающие к $I$, расположены как в стандартном разбиении.
Впрочем, нам это не требуется.)
\end{proof}

\begin{prop}\label{prop:1sharp}
Пусть $I,I'$ --- пара рёбер $\delta$, имеющие общий угол $>\pi$ и другой угол каждого из которых $<\pi$.
Тогда хотя бы один из треугольников разбиения $P$, примыкающий к $I$ или к $I'$, совпадает с треугольником из стандартного разбиения.
\end{prop}

\begin{proof}
Пусть $A$ --- общая вершина $I$ и $I'$.
Тогда один из треугольников разбиения $P$, примыкающих к $I$ и $I'$ и содержащих точку $A$, имеет в $A$ вершину.

Пусть для определённости такой треугольник примыкает к $I$ стороной $b$.
Тогда нам остаётся повторить рассуждение из доказательства предложения~\ref{prop:2sharp} для ребра $I$.
В конце мы получим, что треугольник разбиения $P$, имеющий вершину в другой вершине $I$, расположен как в стандартном разбиении.
\end{proof}

\begin{proof}[Доказательство теоремы~\ref{th:unique}]
Возьмём треугольник $\Delta\in M\setminus\Sigma$.
Предположим, $\Delta$ допускает разбиение $P$, отличающееся от стандартного.
По лемме~\ref{l:n2-cong}, $P$ --- разбиение на $n^2$ треугольников, подобных $\Delta$ с коэффициентом $\frac1n$.

Определим $U$ как объединение треугольников $P$, которые расположены не как в стандартном разбиении.
Тогда, применяя лемму~\ref{l:2or1sharp} и далее предложение~\ref{prop:2sharp} или предложение~\ref{prop:1sharp},
мы видим, что один из треугольников $P$, принадлежащий $U$, совпадает с треугольником из стандартного разбиения.
Но это противоречит построению $U$. 
\end{proof}

\end{document}